\documentstyle[12pt]{article}

\advance\hoffset by -7mm  
\makeatletter
\@addtoreset{equation}{section}
\makeatother

\input amssym.def  
\input amssym.tex

\renewcommand{\title}[1]{\null\vspace{25mm}

\noindent{\Large{\bf #1}}\vspace{10mm}

\noindent {\large By }}
\newcommand{\authors}[1]{\noindent{\large #1}\vspace{3mm}

}
\newcommand{\address}[1]{\noindent #1\vspace{5mm}

}
\renewcommand{\abstract}[1]{\vspace{19mm}

\noindent{\small{\em Abstract.} #1}\vspace{2mm}

}

\begin{document}

\title{Shock formation for forced Burgers equation and application}
\authors{M L Bialy}
\address{School of Math Sciences, Tel-Aviv University, Tel-Aviv 69978,
  Israel\footnote{The author was partially supported by EPSRC. 
    }\\
Email:  bialy@math.tau.ac.il}

\abstract{
We study the inviscid Burgers equation in the presence of spatially
periodic potential force.  We prove that for foliated initial value
problem there are always solutions developing shocks in a finite time.
 We give an application of this result to a quasi-linear system of
 conservation laws which appeared in the study of integrable
 Hamiltonian systems with 1.5 degrees of freedom.
}

\section{Introduction}

Usually it is not a simple task to prove the shock formation for
quasi-linear equations.  We refer to \cite{[Ma]} for an exposition of
techniques.  For the forced Burgers equation there are usually some initial
data which do not lead to formation of shocks.  In this paper we use
integral geometry approach to prove the shock formation for at least
some leaves of foliation which is composed by the graphs of 
solutions of the forced Burgers
equation.
This approach is inspired by Hopf's famouse theorem \cite{[H]} 
in Riemannian geometry,
which can be naturally interpreted as a result on formation of shocks.

\noindent{\large \bf Acknowledgements}

\noindent
The paper was written during my stay at the University of Cambridge.
I would like to thank Robert MacKay for his hospitality and for 
useful discussions.
I would also like to thank Yakov Sinai for helpful discussions on
his paper \cite{[Si]} on Burgers equations.

\section{The main result}

Consider the inviscid Burgers equation
\begin{equation}
f_t + f f_q + F = 0
\label{11}
\end{equation}
We shall assume that the force $F(q,t)=u_q(q,t)$, 
where the potential function $u(q,t)$
satisfies the following requirements.

(2.2a)\quad $u$ is of class $C^2$, 1-periodic in $q, u(q+1, t) =
  u(q,t)$.

(2.2b)\quad for a positive constant $K$, $\int^1_0 u^2_q (q,t) dq<K$
  holds for all $t$.

\noindent
We consider the initial values for $f(q,t)$ depending on parameter $\alpha$
$$
f_\alpha (q,t) |_{t=0} = \varphi (\alpha, q)
$$
where $\varphi$ is a $C^2$-function satisfying:

(2.3a)\quad $\varphi (\alpha, q)$ is a 1-periodic in $q$, for any
  $\alpha$

(2.3b)\quad for any fixed $q$, the mapping $\alpha \mapsto \varphi
  (\alpha, q)$ is a $C^2$-diffeomorphism of $\Bbb R$.

\noindent
Geometrical meaning of (2.3) is that the graphs of $\varphi (\alpha,
q)$ for a $C^2$-foliation of the cylinder ${\Bbb S}^1 \times \Bbb R$.
That is why we refer to such initial data as foliated initial data.

\newtheorem{theorem}{Theorem} 
\begin{theorem}
\label{T1}
Let $u(q,t)$ be the potential of a non-zero force $F$ satisfying
(2.2a,b). Then for any foliated initial data $\varphi (\alpha, q)$ satisfying
(2.3a,b) there always exist the values of $\alpha$ such that the
corresponding solutions of Burgers equation (2.1) develop shock
singularities in a finite (positive or negative) time. 
\end{theorem}

\newtheorem{remark}{Remark}
\begin{remark}
\label{R1}
  The only case when shocks are not created is the case of
  zero force, with the initial data $\varphi (\alpha, q) = \varphi
  (\alpha)$.
 It should be mentioned that
 there are many potentials satisfying (2.2a,b) such that some
  solutions of (2.1) do not form shocks.  For instance, this is always
  the case if $u$ is smooth enough and periodic in both $q$ and $t$.
  In this case KAM theory applies and yields that there are many
  solutions for (2.1) periodic in $q$ and $t$.
\end{remark}

The next result shows that it is not necessarily true that the shocks
described by  Theorem \ref{T1} always appear in a forward time.

\begin{theorem}
\label{T2}
Let $u$ be any $C^2$-function periodic in $q$ satisfying the
following

(2.4a)\quad  $u(q,t) \equiv 0$ for $0 < T \leq t$

(2.4b)\quad for all $0 \leq t \leq T$, $| u_{qq}(q,t)| < \left({\pi
      \over T}\right)^2$

Then there exists a foliated initial data at $t=0$ for (2.1)
such that all the shocks are created in a backward time.
\end{theorem}

\section{Proofs}
For the proof of Theorem \ref{T1} we will need the following
\newtheorem{lemma}{Lemma} 
\begin{lemma}
\label{L1}
Let $V$ be a $C^1$-vector field on $\Bbb R^2$ with the
following property
\begin{equation}
div V  \geq C || V ||^2
\label{21}
\end{equation}
for a positive constant $C$.  Then $V\equiv 0$.
\end{lemma}
\begin{remark}
\label{R2}
This lemma does not generalise to higher dimensions.  There are smooth
vector fields on ${\Bbb R}^n$ for $n \geq 3$ satisfying (3.1)
everywhere.
\end{remark}

\noindent{\bf Proof of Lemma 1}
Integrate (3.1) over the circle $S_r$.  We obtain
\begin{equation}
\int_{S_r} div V d \omega_r \geq C \int_{S_r} ||V 
||^2 d \omega_r
\label{22}
\end{equation}
where $d \omega_r$ is the standard measure on $S_r$.  The left hand
side of (\ref{22}) can be  easily written in the form
\begin{equation}
\int_{S_r} div V d \omega_r = {d \over dr} \int_{B_r}
div V d(vol) = {d \over dr} \left( \int_{S_r} <
  V, n > d \omega_r \right)
\label{23}
\end{equation}
where $\partial B_r = S_r$ and $n$ is a unite normal to
$S_r$.
The right hand side of (\ref{22}) can be estimated by Cauchy--Shwarz
inequality
\begin{equation}
\int_{S_r}||V ||^2 d \omega_r \geq \left(\int_{S_r}
  <V, n>d \omega_r \right)^2 /
\int_{S_r}||n ||^2 d \omega_r = {1 \over 2 \pi r}
  \left(\int_{S_r} <V, n>d \omega_r \right)^2
\label{24}
\end{equation}
Combining (\ref{22}) with (\ref{23}) and (\ref{24}) we obtain the
following
\begin{equation}
\varphi'(r) \geq {C \over 2 \pi r} \quad {\mbox where} \quad
\varphi(r) = \int_{S_r} < V, n> d \omega_r
\label{25}
\end{equation}

It is easy to see that the only solution of (\ref{25}) which is finite
for all $r>0$ is $\varphi \equiv 0$.  But then $V \equiv
0$, by 
(\ref{23}) and (\ref{22}). $\Box$

\noindent
{\bf Proof of Theorem 1}
Proof goes by contradiction.  Let $\varphi (\alpha, q)$ be a foliated
initial data for (\ref{11}) which does not leed to formation of shocks.  
Note that the
characteristics of (\ref{11}) are given by the Newton Equations
\begin{equation}
\left\{ \begin{array}{ll}
\dot{q} & = p \\
\dot{p} & = -u_q
\end{array}
\right.
\label{26}
\end{equation}
The periodicity assumption (2.2a) implies that the flow $g^t$ of (\ref{26})
is complete.  Then we have that the graphs $\left\{ p = f_\alpha (q,t)
\right\}$ form a $C^2$-foliation of the space ${\Bbb R} (p) \times
{\Bbb S}^1 (q) \times {\Bbb R} (t)$.
Define the function $\omega (p,q,t)$ by the rule
$$
\omega \left( f_\alpha (q,t), q,t \right) = {\partial f_\alpha \over
  \partial q} (q,t).
$$
Then $\omega$ is $C^1$, and it is easy to see that the following
equation holds true:
\begin{equation}
\omega_t + p \omega_q - u_q \omega_p + \omega^2 + u_{qq} = 0
\label{27}
\end{equation}
Integrate (\ref{27}) with respect to $q$ over ${\Bbb S}^1$ and obtain
\begin{equation}
- {\partial \over \partial t} \int \omega dq + {\partial \over
  \partial p} \int \omega  u_q dq  = \int \omega^2 dq
\label{28}
\end{equation}

Denote by $V_1 (p,t) = -\int \omega dq, \quad V_2 (p,t) = \int
\omega u_q dq$ and set $ V = (V_1, V_2).$  Then the equation
(\ref{28}) implies for the field $V$ on the plane ${\Bbb R}^2 (p,t)$
satisfies
\begin{equation}
div V = \int \omega^2 dq
\label{29}
\end{equation}
Cauchy-Shwarz inequality applied to the right handside of
 (\ref{29}) together with the
assumption (2.2b) imply that the field $V$ meets the assumption
(\ref{21}) of the lemma.  But then $V$ vanishes
identically and so does $\omega$ (by (\ref{29})) and also $u_q$ (by
(\ref{27})).   This completes the proof. $\Box$

\noindent
{\bf Proof of Theorem 2}
Fix a number $\alpha$ and consider the family 
$M_\alpha$ consisting of those characteristics which for $t \geq T$
can be written
\begin{equation}
q(t, \beta) = \alpha (t-T) + \beta
\label{210}
\end{equation}
It follows that $M_\alpha$ is ordered and defines a smooth foliation of the
semi-cylinder ${\Bbb S}^1 \times \left\{ t \geq 0 \right\}$.  Indeed, the
Jacobi field $\xi_\beta (t) = {\partial q \over \partial \beta} (t,
\beta)$ satisfies the linearised equation
\begin{equation}
\ddot\xi_\beta + u''_{qq} \left( q(t, \beta),t \right) \xi_\beta = 0
\label{211}
\end{equation}
with $\dot\xi_\beta (T)= 0$ by (\ref{210}) 
.  Comparing the equation (\ref{211})
with $\ddot\xi + \left( {\pi \over T}\right)^2 \xi = 0$ on the
segment $[0,T]$ one immediately concludes by (2.4b)
 that $\xi_\beta (t) \not= 0$,
for all t in $[0,T]$.  This implies that ${\partial q \over
  \partial \beta} (t, \beta) > \mbox{const} > 0$.  And thus $M_\alpha$
is a smooth foliation.

To each $M_\alpha$ naturally corresponds the solution $f_\alpha (q,t)$ of
(\ref{11}) defined by the rule
$$
f_\alpha (q(t, \beta),t) = {\partial q \over \partial t} (t, \beta).
$$
The family of solutions $f_\alpha (q,t)$ define the foliated initial
data $\varphi_\alpha (q) = f_\alpha (q,0)$ for which, obviously, the
solutions exist infinite positive time.  The proof is
completed. $\Box$

\section{An Application}

Consider the quasi-linear system $U_t = A(U)U_q$ for $U = (u_1 (q,t)
\dots u_n (q,t))$ of the following form:
\begin{equation}
(u_k)_t = (n-k+1) u_{k-1} (u_1)_q - (u_{k+1})_q \quad for \quad k=1, \dots , n
\label{31}
\end{equation}
where $u_{n+1} \equiv 0$ and $u_0$ is a constant parameter.  This
system naturally appears in the study of polynomial integrals for the
Hamiltonian system (\ref{26}).  It turns out \cite{[B1]} that
(\ref{31}) has a remarkable Hamiltonian form and infinitely many
conservation laws.  It is an open question if all {\it smooth}
solutions of (\ref{31}) defined on the whole cylinder $\Bbb S^1 (q)
\times \Bbb R (t)$ are simple waves solutions (see \cite{[B2]} for
partial results).

We apply Theorem 1 for the elliptic regimes of (\ref{31}), i.e. for
those solutions for which $A(U)$ has no real eigenvalues ($n$ is
automatically even in this case).\footnote{We shall report on strictly
  hyperbolic case elsewhere.}  Let me start with the

\noindent
{\bf Example}
In the case $n=2$, the system has the form
\begin{equation}
\left(
\begin{array}{l}
u_1 \\
u_2 \\
\end{array}
\right)_t =
\left(
\begin{array}{ll}
2u_0 & -1 \\
 u_1  &  0  \\
\end{array}
\right)
\left(
\begin{array}{l}
u_1  \\
u_2  \\
\end{array}
\right)_q
\label{32}
\end{equation}
For a solution $U= (u_1, u_2)$ lying in the elliptic domain,
ie. satisfying $u_1 (q,t) > u^2_0$ introduce the function 
\begin{equation}
E(t) = \int_{S^1} (u_1)^\gamma dq \quad
for \quad \gamma \in (0,1).
\end{equation}
Then the direct computation using (\ref{32}) yields 
$$\ddot E(t) =
\gamma (\gamma -1) \int_{S^1} (u_1)^{\gamma - 2} \left((u_1)^2_t - 2u_0
(u_1)_t (u_1)q + (u_1)^2_q \right)dq $$
Since $u_1>u^2_0$ the integrand is non-negative
and hence by the choice of $\gamma$ in $(0,1)$, one obtains that
the function $E(t)$ is a
concave positive function.  Thus $E \equiv 0$ and then $u_1, u_2$ are
constants.

For higher $n$ we prove the following:

\begin{theorem}
\label{T3}
Let $n$ be even and greater than two and let $U=(u_1 (q,t) \dots u_n (q,t))$ be 
a smooth solution for (\ref{31}) defined on $\Bbb S^1 (q) \times \Bbb
R (t)$ satisfying

(4.3a)\quad $U$ is such that the matrix $A(U)$ has no real eigenvalues. 

(4.3b)\quad $\int_{S^1} u^2_1 (q,t) dq < K$, for some positive
  constant $K$.

\noindent
Then $U \equiv const$.
\end{theorem}

\noindent
{\bf Proof of theorem 3}
For a solution $U$ introduce a polynomial function
$$
F={1 \over n+1} p^{n+1} + u_0 p^n + u_1 p^{n-1} + \cdots + u_n
$$

It can be easily checked that the system (\ref{31}) 
expresses the fact that the function $F(p,q,t)$ satisfies the
equation
$$
F_t + pF_q - (u_1)_q F_p = 0
$$
ie. $F$ has constant values along the Hamiltonian flow of (\ref{26})
(with $u=u_1$). Moreover, it turns out that characteristic polynomial of $A(U)$
satisfies
$$
det (A(U) - \lambda I) = F_p (-\lambda, u_0, u_1, \dots u_n)
$$
But then by the assumption (4.3a) $F_p$ does not vanish and hence the
levels of $F$ determine the foliation consisting of graphs of
solutions of (\ref{11}) with the potential $u_1 (q,t)$.  Then it follows from
Theorem 1 that $U$ must be a constant. $\Box$

\section{Open questions}
We formulate here some natural open questions

1. It is not clear if the growth condition (2.2b) is really essential for
theorem 1 to be true.  Moreover, since theorem 1 is used for the proof
of theorem 3 we had to assume (4.3b).  But the argument in the Example
indicates that probably this assumption may be omitted.

2. The lemma used for the proof of theorem 1 does not generalise to
highest dimensions.  It would be interesting to find some other integral
geometric tools applicable in higher dimensions.  One such tool was
suggested in \cite{[B-P]} for potentials periodic both in space and in
time.

3. It is important to understand if there exist non-zero potential
forces satisfying (2.2 a, b) such that all orbits have no
conjugate points.  Our method does not imply that the force $F$ must vanish
in this case,
though it is close to that.  Such a dichotomy is well known in this
type of question (see for example \cite{[C-K]}).

\hfill

\end{document}